\documentclass[11pt,a4paper]{article}

\usepackage{amsmath}
\usepackage{amssymb,amsfonts,a4}
\usepackage{eucal,bbm,stmaryrd}

\begin{document}

\numberwithin{equation}{section}

\newcommand{\Z}{\mathbbm{Z}}
\newcommand{\R}{\mathbbm{R}}
\newcommand{\N}{\mathbbm{N}}
\newcommand{\Ew}{\mathbbm{E}}
\newcommand{\Bew}{\noindent{\bf Proof:}\;}
\newcommand{\qed}{\hspace*{3em} \hfill{$\square$}}
\newcommand{\weg}{\setminus}
\newcommand{\gdw}{\Longleftrightarrow}
\newcommand{\sgdw}{\Leftrightarrow}
\newcommand{\nach}{\longrightarrow}
\newcommand{\alle}{\;\forall\;}
\newcommand{\gibt}{\;\exists\;}
\newtheorem{satz}{Theorem}
\newtheorem{lemma}{Lemma}

\newcommand{\Bed}{\;|\;}
\newcommand{\F}{\mathcal{F}}
\newcommand{\G}{\mathcal{G}}
\newcommand{\p}{\mathcal{P}}
\newcommand{\E}{\mathcal{E}}
\newcommand{\A}{\mathcal{A}}
\newcommand{\C}{\mathcal{C}}
\newcommand{\V}{\mathcal{L}}
\newcommand{\w}{\mathcal{V}}
\newcommand{\W}{\mathcal{W}}
\newcommand{\Borel}{\mathcal{B}}
\newcommand{\leb}{\lambda^2}
\newcommand{\K}{\mathcal{K}}

\newcommand{\alg}{\mathcal{A}}
\newcommand{\Y}{\mathcal{Y}}
\newcommand{\X}{\mathcal{X}}
\newcommand{\Xd}{\mathcal{X}_{\delta}}
\newcommand{\Xn}{\mathcal{X}_{n}}
\newcommand{\Xrn}{\mathcal{X}_{R}}
\newcommand{\Xbrn}{\mathcal{X}_{R,n}}
\newcommand{\Ye}{\mathcal{Y}}
\newcommand{\Arg}{Y_{\Lambda_n}\bar{Y}_{\Lambda_n^c}}
\newcommand{\xyn}{X\bar{Y}_{\Lambda_n^c}^o}
\newcommand{\ync}{\bar{Y}_{\Lambda_n^c}}
\newcommand{\ybnc}{\bar{Y}}
\newcommand{\yks}{Y{\sigma'}}
\newcommand{\yk}{Y{\sigma}}
\newcommand{\gmu}{\gamma}

%%%%%%%%%%%%%%%%%%%%%%%%%%%%%%%%%%%%%%%%%%%%%%%%%%%%%%%%%%%%%%%%%%%%%%%%

\begin{center}
{\bf \LARGE Two-dimensional Gibbsian point processes 
with continuous spin-symmetries}\\
\vspace{ 1 cm}
Thomas Richthammer\\
Mathematisches Institut der Universit\"at M\"unchen\\
Theresienstra\ss e 39, D-80333 M\"unchen\\
Email: Thomas.Richthammer@mathematik.uni-muenchen.de\\
Tel: +49 89 2180 4633\\
Fax: +49 89 2180 4032 
\end{center}

\vspace{ 1cm} 

\begin{abstract}
We consider two-dimensional marked point processes which are Gibbsian
with a two-body-potential of the form $U = JV + K$, where $J$ and $K$
depend on the positions and $V$ depends on the marks of the two 
particles considered. $V$ is supposed to have a continuous symmetry. 
We will generalise the famous Mermin-Wagner-Dobrushin-Shlosman 
theorem to this setting in order to show that the Gibbsian process 
is invariant under the given symmetry, when instead of smoothness 
conditions only continuity conditions are assumed. 
We will achieve this by using
Ruelle's superstability estimates and percolation arguments. 
\end{abstract}

\section{Introduction}

\begin{sloppypar}

Gibbsian processes were introduced by  R.~L.~Dobrushin (see \cite{D68}
and \cite{D70}), O.~E.~Lanford and D.~Ruelle (see \cite{LR}) as a
model for equilibrium states in statistical physics. (For general
results on Gibbs measures on a d-dimensional lattice we refer to the
detailed book of H.-O.~Georgii \cite{G}, which covers a wide range of
phenomena.) The first results concerned existence and uniqueness of
Gibbs measures and the structure of the set of Gibbs measures related
to a given potential. The question of uniqueness is of special
importance, as the nonuniqueness of Gibbs measures can be interpreted
as a certain type of phase transition occurring within the particle
system. A phase transition occurs whenever a symmetry of the
potential is broken, so it is natural to ask, under which
conditions symmetries are broken or conserved. The answer to this
question depends on the type of the symmetry (discrete or continuous),
the number of spatial dimensions and smoothness and decay conditions
on the potential (see \cite{G}, chapters 6.2, 8, 9 and 20). 
It turns out that the case of continuous
symmetries in two dimensions is especially interesting. The first
progress in this case was achieved by M.~D.~Mermin and H.~Wagner, who
showed for special two-dimensional lattice models that symmetries are
conserved (\cite{MW} and \cite{M}).
In \cite{DS} R.~L.~Dobrushin and S.~B.~Shlosman established 
conservation of symmetries for more general 
potentials which satisfy smoothness and decay conditions, and
C.-E.~Pfister improved this result in  \cite{P}; considering a
continuous model J.~Fr\"ohlich and C.~Pfister (\cite{FP}) worked 
with the concept of superstability (see \cite{R}), 
whereas H.-O.~Georgii gave a fairly elementary proof in \cite{G99}. 
All these results rely on the smoothness of the interaction, and
only recently D.~Ioffe, S.~Shlosman and Y.~Velenik showed that
mere continuity suffices in the lattice model (\cite{ISV}) using
a perturbation expansion and  percolation theory.\\  

We will generalise the last result from a lattice to a continuous
model, using superstability techniques. Apart from that we will mimic 
the proof of D.~Ioffe, S.~Shlosman and Y.~Velenik 
and use a very similar percolation argument.\\

In section 2 we will describe the situation considered 
and state the result obtained. The precise setting is then given 
in section 3. In section 4 a proof of a weaker version of the result 
is given. The proofs of all lemmas are relegated to section 5, and in
section 6 we will show how to deal with the general case.\\

Acknowledgement: I would like to thank Prof.~Dr.~H.-O. Georgii 
for suggesting the problem and many helpful comments.

%%%%%%%%%%%%%%%%%%%%%%%%%%%%%%%%%%%%%%%%%%%%%%%%%%%%%%%%%%%%%%%%%%%

\section{The Result}

We consider infinitely many particles in the plane, 
where a particle has a position in $\R^2$ and internal degrees
of freedom. These can be modeled by assigning to the particle a value from
some measurable spin space (or mark space) $S$. The particles may interact via 
a pair potential $U$. So  $U$ is a measurable function 
\begin{displaymath}
U : (\R^2 \times S)^2 \; \to \; \bar{\R} \; := \; \R \cup \{+\infty\}
\end{displaymath}
such that $U(y_1,y_2) = U(y_2,y_1)$ for all $(y_1,y_2) \in D$,
i. e. $U$ is symmetric. Here we assume $U$ to be of the form 
\begin{equation} \label{psidompotsup}
U(x_1,\sigma_1;x_2,\sigma_2) \; = \; J(x_1-x_2)\tilde{U}(\sigma_1,\sigma_2)
+ K(x_1 - x_2) 
\end{equation}
such that the functions $\tilde{U}: S^2 \to \R$, $J: \R^2 \to \R$ and 
$K: \R^2 \to \bar{\R}$ are measurable and symmetric, and $J$ is
$\psi$-dominated, i. e. 
\begin{displaymath}
|J(x)| \;  (1 + \|x\|^2) \;\le\; \psi(\|x\|)  \qquad \alle x \in \R^2, 
\end{displaymath}
where $\psi: \R_+ := [0,\infty[ \to \R_+$ is a given decreasing function
such that
\begin{equation*}
\int_0^{\infty} \psi(r)r \; dr\; := \; \psi_s  < \infty.
\end{equation*}
We will call a potential $U$ of the above form \eqref{psidompotsup} 
a $\psi$-dominated potential corresponding to $J,K,\tilde{U}$.\\ 
   
We are only interested in the equilibrium states of a thermodynamical system 
as described above, and as a model for these we take the concept of Gibbs
measures. Supposing that the given potential has some internal symmetry,
we would like to know whether the possible equilibrium states inherit this
symmetry necessarily. 
For example, considering a potential which does not change under
rotation of spins, under what conditions are the equilibrium states invariant
under spin rotation? Here we are concerned with continuous symmetries only,
so that we can model the symmetries by a Lie-group $G$ acting on the spin
space $S$. Our result is then the following:    

\begin{satz} \label{symallgsup}
Let $(S,\Borel(S),\lambda_S)$ be a probability space such that $S$ is
a compact topological space and $\Borel(S)$ its
Borel-$\sigma$-algebra.  
Let $G$ be a compact connected Lie-group operating on $S$ such that
the operation is continuous and the reference measure $\lambda_S$
is $G$-invariant. Let $U$ be a superstable, lower regular,
$\psi$-dominated potential corresponding to $J,K,\tilde{U}$ 
such that  $\tilde{U}$ is continuous and $G$-invariant.
Then every tempered Gibbs measure corresponding to $U$ is
$G$-invariant.
\end{satz}
The exact definitions of the objects and properties in the formulation 
of the above theorem will be given in the next section.

%%%%%%%%%%%%%%%%%%%%%%%%%%%%%%%%%%%%%%%%%%%%%%%%%%%%%%%%%%%%%%%%%%%%%%
%%%%%%%%%%%%%%%%%%%%%%%%%%%%%%%%%%%%%%%%%%%%%%%%%%%%%%%%%%%%%%%%%%%%%%

\section{The Setting}

\subsection{Configurations of particles}

We consider the plane $\R^2$ with maximum norm $\|.\|$. Let 
\begin{displaymath}
 \Lambda_t \; := \; [-t,t[^2 \quad \text{ for } t \in \R_+ \quad 
\text{ and } \qquad  
C_r \; := \; r + \left[ -\frac{1}{2},\frac{1}{2} \right[^2  \quad
 \text{ for } \quad r \in \Z^2
\end{displaymath}
be subsets of $\R^2$. 
On $\R^2$ let $\Borel^2$ be the Borel-$\sigma$-algebra, and 
$\Borel^2_b \subset \Borel^2$ the set of all bounded Borel sets.
The Lebesgue measure on $(\R^2,\Borel^2)$ will be denoted by $\leb$.\\

For describing the marks or spins of the particles let $S$ be a
topological space, $\Borel(S)$ the Borel-$\sigma$-algebra on $S$ and
$\lambda_S$ a normed reference measure on $(S,\Borel(S))$. As
$\lambda_S$ is the only measure to be considered on $(S,\Borel(S))$, we 
will simply write  $d\sigma := d\lambda_S(\sigma)$ 
when integrating with respect to $\lambda_S$.\\ 

A configuration $Y$ of marked particles is described by a subset of 
$\R^2 \times S$ which is a locally finite, in that 
$|Y \cap (\Lambda \times S)| < \infty$ for all $\Lambda \in \Borel^2_b$, 
and simple, in that for all $ (x_1,\sigma_1) \ne (x_2,\sigma_2) \in Y$
we have $x_1 \ne x_2$. The configuration space $\Y$ is defined to be 
the set of all locally finite and simple subsets of $\R^2 \times S$.  
A configuration $Y \in \Y$ is said to be finite if $|Y| < \infty$.
Given a particle $y \in \R^2 \times S$, we want to consider the
position $y^o \in \R^2$ and the spin $\sigma_y \in S$ of the particle,
and given a configuration $Y \in \Y$ let $Y^o := \{ x \in \R^2: \gibt
\sigma \in E: (x,\sigma) \in Y\}$. For $Y \in \Ye$ and $x \in \R^2$  
such that $(x,\sigma) \in Y$ let $\sigma_x(Y) := \sigma$ and 
$\sigma_x := \sigma_x(Y)$ if it is
clear which configuration is to be considered.\\

For $Y \in \Y, \Lambda \in \Borel^2, B \in \Borel(S)$ let $Y_{\Lambda,B}
:= Y \cap (\Lambda \times B)$ and $Y_{\Lambda} := Y_{\Lambda,S}$ the
restriction of $Y$ to $\Lambda \times S$ and $\Lambda$ respectively,  
$\Y_{\Lambda} := \{ Y \in \Y: Y \subset \Lambda \times S\}$ the set of
all configurations in $\Lambda$, $N_{\Lambda,B}(Y) :=
|Y_{\Lambda, B}|$ the number of particles of $Y$ in $\Lambda$ with marks
in $B$ and $N_{\Lambda} := N_{\Lambda,S}$. 
The counting variables $N_{\Lambda, B}$ generate a
$\sigma-$algebra $\F_{\Y}$ on $\Y$.  
For $\Lambda \in \Borel^2$ let $\F'_{\Y,\Lambda}$ be the
$\sigma$-algebra on $\Y_{\Lambda}$ obtained by restricting
$\F_{\Y}$  to $\Y_{\Lambda}$, and let 
$\F_{\Y,\Lambda} := e_{\Lambda}^{-1} \F'_{\Y,\Lambda}$ be the 
$\sigma$-algebra on $\Y$ obtained from $\F'_{\Y,\Lambda}$
 by the restriction mapping 
$e_{\Lambda}: \Y \to \Y_{\Lambda}, Y \mapsto Y_{\Lambda}$. 
For disjoint sets 
$\Lambda_1,\Lambda_2 \in \Borel^2$ and configurations 
$Y,\bar{Y} \in \Y$ let 
$Y_{\Lambda_1}\bar{Y}_{\Lambda_2} := Y_{\Lambda_1} \cup
\bar{Y}_{\Lambda_2}$ .
\\

The mean quadratic
particle density per unit square for $Y \in \Y$ is defined by  
\begin{equation*} \label{dichtnsup}
s_n(Y) \; := \; \frac{1}{\leb(\Lambda_{n+\frac{1}{2}})} 
\sum_{r \in \Z^2 \cap \Lambda_{n+\frac{1}{2}}} N_{C_r}^2(Y). 
\end{equation*}
A configuration $Y \in \Y$ is said to be tempered if
$s(Y) := \sup_{n \in \N} s_n(Y) < \infty$.
Let $\Y_t \in \F_{\Y}$ be the set of all tempered configurations.\\

Now similar objects can be considered for particles without marks. Let
$\X \; := \; \{ X \subset \R^2: \;  
|X \cap \Lambda | < \infty \alle \Lambda \in \Borel^2_b \}$
be the configuration space of particle positions.  
The restrictions $X_{\Lambda}$, the set of configurations in
$\Lambda\; $ $\X_{\Lambda}$, the counting variables
$N_{\Lambda}$, the $\sigma-$algebras $\F_{\X}$, $\F'_{\X,\Lambda}$ and  
$\F_{\X,\Lambda}$ and $X_{\Lambda_1}\bar{X}_{\Lambda_2}$ are then
defined analogously to the objects above. 
The projection $o: \Y \to \X$, $Y \mapsto Y^o$ obviously is measurable, so 
$\F_{\X}$ can be considered as a subset
of $\F_{\Y}$ via the identification of a set $\X_1 \in \F_{\X}$ with 
$o^{-1} \X_1 \in \F_{\Y}$. 
For example we have that $\Y_t \in \F_{\X}$.
For any $X \in \X$ and a family of marks $(\sigma_x)_{x \in X}$ let
$(X,\sigma) := \{ (x,\sigma_x): x \in X \}$ the configuration determined
by $X$ and $\sigma$.\\

Let $z > 0$ be an activity parameter which will be fixed 
throughout this paper. Let $\nu := \nu_z$ be the distribution 
of the Poisson point process on $(\Y,\F_{\Y})$ with intensity $z$
and distribution of marks $\lambda_S$, and $\nu^o := \nu^o_z$ be the
distribution of the Poisson point process on $(\X,\F_{\X})$ with 
intensity $z$. So 
\begin{equation*} 
\int f d\nu^o \quad = \quad e^{-z \leb(\Lambda)} \; \sum_{ k \ge 0} \;  
 \frac{z^k}{k!} \; \int_{\Lambda^k} dx_1 ... dx_k \; 
 f(\{x_i : 1 \le i  \le k\} ), 
\end{equation*}
for any $\F_{\X,\Lambda}-$measurable nonnegative function $f: \X \to
\R_+$ and 
\begin{equation*} 
\int f d\nu \quad = \quad \int \nu^o(dX) \int_{S^{X_{\Lambda}}} 
d\sigma_{X_{\Lambda}} \; f( (X_{\Lambda},\sigma) ), 
\end{equation*}
for any $\F_{\Y,\Lambda}-$measurable nonnegative function 
$f: \Y \to \R_+$. 

%%%%%%%%%%%%%%%%%%%%%%%%%%%%%%%%%%%%%%%%%%%%%%%%%%%%%%%%%%%

\subsection{Configurations of bonds}

For any set $Z$ and distinct $z_1,z_2 \in Z$ let 
$z_1z_2 :=\{ z_1,z_2\}$ be the bond joining $z_1$ and $z_2$. Let 
$E(Z) \; := \; \{ z_1z_2 : \; z_1,z_2 \in Z, z_1 \ne z_2 \}$
be the set of all bonds in $Z$.  On $E(\R^2)$ the $\sigma$-algebra
\begin{displaymath}
\F_{E(\R^2)} \; := \; \{ \{ x_1x_2 \in E(\R^2):  
(x_1,x_2) \in B \} :  \; B \in (\Borel^2)^2 \}
\end{displaymath} 
is given. Let 
\begin{displaymath}
\E \; := \;  \{ E \subset E(\R^2): \; |\{ xy \in E : 
xy \subset B \} | < \infty \; \alle B \in \Borel_b^2  \}
\end{displaymath}
be the configuration space of bonds, i. e. the set of all locally 
finite bond sets. On $\E$ the $\sigma$-algebra 
$\F_{\E}$ is defined to be generated by the counting variables 
$N_{E'}: \E \to \N,  E \mapsto | E' \cap E|$ 
$\;(E' \in \F_{E(\R^2)})$.\\   

For a countable set $E \in \E$ one can also consider the 
Bernoulli-$\sigma$-algebra $\Borel_E$ on $\E_{E} := \p(E) \subset \E$,  
which is defined to be generated by the family of sets 
$(\{E' \subset E: e \in E'\})_{e \in E}$. 
It is easy to check that the inclusion 
$(\E_E,\Borel_E) \to (\E, \F_{\E})$ is measurable.   
Thus any probability measure on  
$(\E_E,\Borel_E)$ can trivially be extended to $(\E, \F_{\E})$.\\

Given a countable set $E$ and a family 
$(\epsilon_e)_{e \in E}$ of real numbers in $[0,1]$ the Bernoulli
measure on $(\p(E),\Borel_E)$ is defined as the unique probability measure  
for which the events $(\; \{E' \subset E : e \in E' \}\;)_{e \in E}$ 
are independent with probabilities $(\epsilon_e)_{e \in E}$.

%%%%%%%%%%%%%%%%%%%%%%%%%%%%%%%%%%%%%%%%%%%%%%%%%%%

\subsection{Interaction and superstability}

Our next step is to introduce the interaction between particles. 
As mentioned before we will consider a $\psi$-dominated 
potential corresponding to $J,K,\tilde{U}$ as defined in and below of
\eqref{psidompotsup}. The energy of a finite configuration $Y \in \Y$ is
defined as
\begin{equation*} %\label{energiesup}
H^U (Y) \; := \sum_{y_1y_2 \in E(Y)} U(y_1,y_2), 
\end{equation*}
and for two finite configurations $Y,Y' \in \Y$ let 
\begin{equation} \label{wechselsup}
W^U (Y,Y') \; := \sum_{y_1 \in Y} \sum_{y_2 \in Y'} U(y_1,y_2)
\end{equation}
be the interaction energy of the configurations. Definition 
\eqref{wechselsup} can be extended to infinite configuration $Y'$ 
whenever $W^U (Y,Y'_ {\Lambda})$ converges   
as $\Lambda \uparrow \R^2$ through the net $\Borel^2_b$.\\

For a configuration $Y \in \Y$ let 
$\Z^2(Y) \; := \; \{ r \in \Z^2: N_{C_r}(Y) > 0 \}$
be the minimal set of lattice points such that the corresponding
squares cover $Y$. Then a potential is called superstable if there
are real constants $A > 0$ and $B \ge 0$  such that for all finite
configurations $Y \in \Y$ 
\begin{equation*} %\label{supstabsup}
H^U (Y) \; 
\ge \;\sum_{r \in \Z^2(Y)} \;[A N_{C_r}(Y)^2 - B N_{C_r}(Y)].
\end{equation*}
A potential is called lower regular if there is a decreasing function
$\Psi: \N \to \R_+$ such that 
\begin{equation*}
\sum_{r \in \Z^2} \Psi(\|r\|) \; < \; \infty \quad
\text{ and }
\end{equation*}
\begin{equation*} %\label{regabsup}
W^U (Y,Y') \; \ge \; -\sum_{r \in \Z^2(Y)} 
\sum_{s \in \Z^2(Y')} \; \Psi(\|r-s\|) \; 
[\frac{1}{2} N_{C_r}(Y)^2 + \frac{1}{2} N_{C_s}(Y')^2] 
\end{equation*}
for all finite configurations $Y,Y' \in \Y$. 
Note that any $\psi$-dominated potential corresponding to 
$J,K,\tilde{U}$ such that also 
$K(x) \ge - \psi(\|x\|)$ for all $x \in \R^2$ is lower regular.\\

It is well known that
for any superstable and lower regular potential $U$, any finite
configuration $Y \in \Y$ and any tempered configuration $Y' \in \Y_t$ 
the interaction energy $W^U(Y,Y')$ exists in $]-\infty, \infty]$, see 
\cite{R} for example.

%%%%%%%%%%%%%%%%%%%%%%%%%%%%%%%%%%%%%%%%%%%%%%%%

\subsection{Gibbs measures}

Given a superstable and lower regular potential $U$, 
the Hamiltonian of a configuration $Y \in \Y$ in $\Lambda \in \Borel^2_b$ 
with boundary condition $\bar{Y}_{\Lambda^c} \in \Y_t$ is defined by 
\begin{equation*} 
H^U_{\Lambda}(Y_{\Lambda}\bar{Y}_{\Lambda^c}) \; := \; 
H^U (Y_{\Lambda}) +  W^U(Y_{\Lambda},\bar{Y}_{\Lambda^c}) \; =  
\sum_{y_1y_2 \in E(Y_{\Lambda}\bar{Y}_{\Lambda^c}): \; 
 y_1^o y_2^o \cap \Lambda \ne \emptyset} U(y_1,y_2).
\end{equation*}
The integral 
\begin{equation*} %\label{zustandsup}
Z^U_{\Lambda}(\bar{Y}) \; := \; \int
\nu(dY) \; e^{-H^U_{\Lambda}(Y_{\Lambda}\bar{Y}_{\Lambda^c})}
\end{equation*}
is called the partition function in $\Lambda \in \Borel^2_b$ 
for the boundary condition $\bar{Y}_{\Lambda^c} \in \Y_t$. 
Using superstability and lower
regularity of $U$ and temperedness of $\bar{Y}$ one can show that 
$Z^U_{\Lambda}(\bar{Y})$ is finite (see \cite{R} for example), 
and considering the empty
configuration $Y$ one can show that $Z^U_{\Lambda}(\bar{Y})$ is positive. 
The Gibbs distribution $\gmu^U_{\Lambda}(.|\bar{Y})$ in
$\Lambda \in \Borel_b$ with boundary condition $\bar{Y}_{\Lambda^c} 
\in \Y_t$, potential $U$ and activity $z$ is thus well defined by 
\begin{displaymath} 
\gmu^U_{\Lambda}(A|\bar{Y}) \; := \; 
Z^U_{\Lambda}(\bar{Y})^{-1} \int \nu(dY) \; 
e^{-H^U_{\Lambda}(Y_{\Lambda}\bar{Y}_{\Lambda^c})}
1_A(Y_{\Lambda}\bar{Y}_{\Lambda^c}) \quad 
\text{ for } \quad  A \in \F_{\Y}.
\end{displaymath} 
$\gmu^U_{\Lambda}$ is a probability kernel from 
$(\Y,\F_{\Y})$ to $(\Y,\F_{\Y})$. Let 
$\gmu_{\Lambda} := \gmu^U_{\Lambda}$ if it is clear which potential is
considered. Let  
\begin{displaymath}
\begin{split}
\G(U) \quad := \quad \{ \mu \in &\p(\Y,\F_{\Y}) : \; \mu(\Y_t) = 1\\
&\mu(A|\F_{\Y,\R^2 \weg \Lambda}) = 
 \gmu^U_{\Lambda}(A|.) \text{ $\mu$-a.s. } \alle A \in \F_{\Y},
\Lambda \in \Borel^2_b \}
\end{split}
\end{displaymath}
be the set of all tempered Gibbs measures for the potential $U$ 
and the activity $z$. It is easy to see 
that for any probability measure $\mu \in \p(\Y,\F_{\Y})$
such that $\mu(\Y_t)=1$ the equivalence 
\begin{equation*} 
\mu \in \G(U) \quad \gdw \quad (\mu \gmu^U_{\Lambda} = \mu
 \alle \Lambda \in \Borel^2_b)
\end{equation*} 
holds. So for every $\mu \in \G(U)$, $f: \Y \to \R_+$ measurable and  
$\Lambda \in \Borel^2_b$ we have
\begin{equation} \label{gibbsaequsup}
\int \mu(dY) \; f(Y)\; = \; 
\int \mu(d\bar{Y}) \int \gmu_{\Lambda}^U (dY | \bar{Y}) \; 
f(Y_{\Lambda}\bar{Y}_{\Lambda^c}). 
\end{equation}

For a superstable and lower regular potential $U$ and a tempered Gibbs 
measure $\mu \in \G(U)$, the correlation function
$\rho^{U,\mu}$ of $\mu$ is defined by 
\begin{equation*} %\label{korsup}
\rho^{U,\mu}(Y) \; = \; e^{-H^U(Y)} \int \mu (d\bar{Y}) \;  
     e^{- W^U (Y,\bar{Y})}
\end{equation*}
for any finite configuration $Y$. It is a remarkable consequence of
Ruelle's superstability estimates that there is a constant $\xi \in \R$
such that  
\begin{equation} \label{korbeschrsup}
\rho^{U,\mu}(Y) \quad \le \quad \xi^{|Y|}
\end{equation}
for any finite configuration $Y \in \Y$. (For a proof see \cite{R}.)
We will call a $\xi \in \R$ satisfying \eqref{korbeschrsup} a
Ruelle bound. Actually we will need this bound on the correlation
function in the following way: 
\begin{lemma} \label{lekorabsup}
Let $U$ be a superstable and lower regular potential, 
$\mu \in \G(U)$ a tempered Gibbs measure and $\xi \in \R$ a
Ruelle bound. Then we have 
\begin{equation} \label{korabsup}
\int \mu(dY) \sideset{}{^{\neq}}\sum_{x_1,...,x_m \in Y^o} 
f(x_1,...,x_m) \quad
\le \quad (z\xi)^m \; \int dx_1 ... dx_m \; f(x_1,...,x_m)
\end{equation}
for every integer $m \ge 0$ and every measurable function
$f: (\R^2)^m \to \R_+$. 
\end{lemma}
We use $\Sigma^{\neq}$ as a shorthand notation for a multiple sum 
such that the summation indices are assumed to be pairwise distinct.     

%%%%%%%%%%%%%%%%%%%%%%%%%%%%%%%%%%%%%%%%%%%%%%%%%%%%%%

\subsection{Transformations of spins}

Now let the spin space $S$ be a compact topological space,  
and $G$ be a compact, connected Lie-group operating on $S$, 
\begin{displaymath}
op: G \times S \to S, \qquad (\tau,\sigma) \mapsto op(\tau,\sigma) =: 
\tau(\sigma),
\end{displaymath}
such that the operation is measurable.\\ 

For every  $\tau \in G$ we also consider   
$\tilde{\tau}: \Y \to \Y$, $\tilde{\tau}(Y) =
\{(x,\tau(\sigma)): (x,\sigma) \in Y\}$,  and 
$\bar{\tau}: D \to D$, $\bar{\tau}(x_1,\sigma_1;x_2,\sigma_2) =
(x_1,\tau(\sigma_1);x_2,\tau(\sigma_2))$. Usually these mappings will
again be denoted by  $\tau$.
Furthermore, for a configuration $Y \in \Y$ and $\tau: Y^o \to G$ 
we write $\tau Y := \tau(Y) := \{ (x,\tau(x) (\sigma)): 
(x,\sigma) \in Y\}$.\\

$\tau \in G$ is called a symmetry of a given pair potential $U$
if $U \circ \tau = U$. If this holds for every $\tau \in G$, then
$U$ is said to be $G$-invariant. The reference measure $\lambda_S$ is
called $G$-invariant if $\lambda_S \circ \tau^{-1} = \lambda_S$ for all  
$\tau \in G$, and a Gibbs measure $\mu \in \G(U)$ is called
$G$-invariant if $\mu \circ \tau^{-1} = \mu$ for all $\tau \in G$.

%%%%%%%%%%%%%%%%%%%%%%%%%%%%%%%%%%%%%%%%%%%%%%%%%%%%%%%%%%%%%
%%%%%%%%%%%%%%%%%%%%%%%%%%%%%%%%%%%%%%%%%%%%%%%%%%%%%%%%%%%%%

\section{The case of $S^1$-action}

We will first consider the mark space 
$(S,\Borel(S),\lambda_S) := (S^1, \Borel(S^1),\lambda_{S^1})$, 
where $S^1$ is the unit circle, $\Borel(S^1)$ is the 
Borel-$\sigma-$algebra on $S^1$ and $\lambda_{S^1}$ 
is the Lebesgue-measure on $S^1$, and transformations 
$\tau \in G := \{\tau_{\sigma}: \sigma \in S^1\}$,  
where $\tau_ {\sigma}$ is defined to be the rotation with angle 
$\sigma$. For $\sigma, \sigma' \in S^1 = \R / (2\pi\Z)$ we write 
$\tau_{\sigma}(\sigma') =: \sigma'+\sigma$.
In order to simplify notation we identify a rotation with
its angle, i. e. we identify $S^1 = \R / 2\pi\Z$ with $[0,2\pi[$, 
and so we consider functions on $S^1$ as $2\pi$-periodic functions on
$\R$ whenever possible.\\ 

If all rotations $\tau \in G$ are symmetries of the
$\psi$-dominated potential $U$ corresponding to $J,K,\tilde{U}$, 
then $U$ can also be written in the form 
\begin{displaymath}
U(x_1,\sigma_1;x_2,\sigma_2) = J(x_1-x_2)V(\sigma_1-\sigma_2)
+ K(x_1 - x_2),
\end{displaymath}
where $V: S \to \R$ is defined by $V(\sigma) := \tilde{U}(\sigma,0)$. 
On the other hand a potential of the above form is $G$-invariant. 
It is called the $\psi$-dominated potential corresponding to $J,K,V$. 
As an additional preliminary simplification we assume
that $J \ge 0$. So we consider the follwing special case of 
theorem \ref{symallgsup}: 
\begin{satz} \label{symsup}
Let $U$ be a superstable, lower regular 
$\psi$-dominated potential corresponding to $J,K,V$ such that 
$J \ge 0$ and $V$ is continuous.\\
Then every tempered Gibbs measure corresponding to $U$ 
is $G$-invariant. 
\end{satz}
In the following subsections we will give a proof of this theorem.

%%%%%%%%%%%%%%%%%%%%%%%%%%%%%%%%%%%%%%%%%%%%%%%%%%%%%%%%%%%%%%%%%%%%%%

\subsection{Constants and Decomposition of $V$}  \label{constantssup}

Let $U$ be a potential with the properties stated in Theorem \ref{symsup}, 
$\mu \in \G(U)$ a tempered Gibbs measure and $\xi \in \R$ a Ruelle bound
satisfying \eqref{korbeschrsup} and $1 < 2z\xi$, 
where again $z$ is the intensity of the underlying Poisson point process. 
As a consequence of the $\psi$-domination of $J$ and the integrability
condition on $\psi$ there is a real constant $c_J$ such that
\begin{equation*}
1 \; + \; \psi(0) \; + \;  \int J(x) (1 + \|x\|^2) dx \quad \le \quad c_J,
\end{equation*} 
and there are real constants $c(R)$ for $R \ge 0$  such that 
$\lim_{R\to \infty} c(R) = 0$ and for all $R \ge 0$ 
\begin{equation} \label{crsup}
\int 1_{\{|x| \ge R\}} \;  J(x) dx \quad \le \quad  c(R). 
\end{equation}
We want to show the $G$-invariance of $\mu$ by an argument similar to 
the one given in \cite{G}, chapter 9.1, proposition (9.1). 
So we fix a transformation $\tau \in ]0,\pi[$, a test cylinder event 
$B \in \F_{\Y,\Lambda_{n'}}$ $(n'\in \N)$ and a real $\delta > 0$.
Furthermore let  $1 > \epsilon > 0$ such that
\begin{equation} \label{bedepssup}
 c_J \epsilon \; < \;  2 c_J z \xi \epsilon \; < \; 1,
\end{equation}
As the above parameters are fixed for the whole proof we will ignore the
dependence of any variable on any of the above parameters.\\

As $V$ is a continuous function on $S^1$, $V$ can be approximated by
trigonometric polynomials due to the Weierstra{\ss} theorem. 
So we have the decomposition $V = \tilde{V} - \tilde{v}$, 
such that $\tilde{V}$ is smooth (i. e. twice continuously
differentiable), and $|\tilde{v}| < \frac{\epsilon}{2}$. 
Defining $v := \tilde{v} + \frac{\epsilon}{2}$ and 
$\bar{V} := \tilde{V} + \frac{\epsilon}{2}$ we get the decomposition
\begin{equation*}
V = \bar{V} - v \qquad \text{ with smooth $\bar{V}$ and } 
0 < v < \epsilon.
\end{equation*}
By symmetrizing $\bar{V}$ and $v$ we can assume $\bar{V}$ and $v$ to
be symmetric. Let $\bar{U}$ be the $\psi$-dominated potential
corresponding to $J,K,\bar{V}$.

%%%%%%%%%%%%%%%%%%%%%%%%%%%%%%%%%%%%%%%%%%%%%%%%%%%%%%%%%%%%%%%%%%%%

\subsection{Decomposition of $\mu$ and the bond process}

For $n \in \N$ and $X \in \X$ we consider the bond set  
\begin{displaymath}
E(X,n) \; := \;  \{ x_1x_2 \in E(X): J(x_1-x_2) \ne 0, 
x_1x_2 \cap \Lambda_n \ne \emptyset \}.
\end{displaymath}
In order to be able to extend the decomposition of the potential
function $V$ to a decomposition of the Hamiltonian we need :
\begin{lemma} \label{umordsup}
For each $n \in \N$ there is a set $\Xn \in \F_{\X}$ such that 
$\mu(\Xn) = 1$ and  
\begin{displaymath}
\sum_{x_1x_2 \in E(X,n)} J(x_1-x_2) \;  < \; \infty
\qquad \alle X \in \Xn. 
\end{displaymath}
\end{lemma}
Now let $n \in \N$ and $Y \in \Xn$ be fixed. 
Because of lemma \ref{umordsup} we have  
\begin{equation*}
H^U_{\Lambda_n}(Y) \quad = \quad  H^{\bar{U}}_{\Lambda_n}(Y) - 
 \sum_{x_1x_2 \in E(Y^o,n)} J(x_1-x_2) \; 
v(\sigma_{x_1}(Y) -  \sigma_{x_2}(Y)),  
\end{equation*}
and therefore
\begin{equation} \label{azerlsup}
e^{-H^U_{\Lambda_n}(Y)} \quad 
 = \sideset{}{'}\sum_{A \subset E(Y^o,n)} \w_n(A,Y), 
\end{equation}
where we have used the shorthand notation 
\begin{equation*}
\w_n(A,Y) \quad := \quad e^{ - H^{\bar{U}}_{\Lambda_n}(Y)} 
 \prod_{x_1x_2 \in A}
 [e^{J(x_1-x_2)\,v(\sigma_{x_1}(Y)-\sigma_{x_2}(Y))}-1]
\end{equation*}
for $n \in \N$, $Y \in \Ye $ and finite $A \subset E(Y^o,n)$.
The summation symbol $\sum'$ in \eqref{azerlsup} indicates that the sum
extends over finite subsets only. For $n \in \N$, $X \in \X_{\Lambda_n}$,
$\bar{Y} \in \Y_{\Lambda_n^c}$ such that $X\bar{Y}^o \in \X_n$, 
finite $A \subset E_n := E(X\bar{Y}^o,n)$, $\E' \in \Borel_{E_n}$ 
and $D \in \F_{\Y}$ we define 
\begin{displaymath}
\begin{split}
\W_n(X,\ybnc) \quad 
 &:= \quad  \int d\sigma_X \; 
 e^{ -H^U_{\Lambda_n}((X,\sigma)\ybnc)} \\
\W_n(A,X,\ybnc) \quad &:= \quad  
 \int d\sigma_X \; \w_n(A,(X,\sigma)\ybnc)\\
\pi_n(\E'|X,\ybnc) \quad &:= \quad \sideset{}{'}\sum_{A \in \E'} 
 \frac{\W_n(A,X,\ybnc)}{\W_n(X,\ybnc)}\\ 
\alpha_n(D|A,X,\ybnc) \quad &:= \quad 
 \frac{1}{\W_n(A,X,\ybnc)} \int d\sigma_X  \; \w_n(A,(X,\sigma)\ybnc) 
\;  1_D ((X,\sigma)\ybnc).  
\end{split}
\end{displaymath}
As $J$ and $v$ are nonnegative the above factors and integrands are
nonnegative, too, and so all products and integrals are well defined.
If $\W_n(X,\ybnc) = 0$ or $X\ybnc^o \notin \X_n$ 
we define $\pi_n(.|X,\ybnc)$ to be the
probability measure on $(\E_{E_n},\Borel_{E_n})$ with whole weight on
the empty set. If  $\W_n(A,X,\ybnc) = 0$ or $X\ybnc^o \notin \X_n$ or 
$A \in \E$ is not a finite
subset of $E_n$ let $\alpha_n(.|A,X,\ybnc)$ be an
arbitrary fixed probability measure on $(\Y,\F_{\Y})$. 
For $n \in \N$, $X \in \X_{\Lambda_n}$ and $\bar{Y} \in \Y_{\Lambda_n^c}$ 
such that $\W_n(X,\ybnc) > 0$ and $X\ybnc^o \in \X_n$ 
we have by \eqref{azerlsup}  
\begin{displaymath}
\pi_n(\E_{E_n} | X,\ybnc) \quad 
=\quad  \frac{1}{\W_n(X,\ybnc)} \; 
 \sideset{}{'}\sum_{A \subset E_n} 
\int d\sigma_X \; \w_n(A,(X,\sigma)\ybnc)  \quad = \quad 1.
\end{displaymath} 
Therefore $\pi_n(.|X,\ybnc)$ is a probability measure on 
$(\E_{E_n},\Borel(\E_{E_n}))$ and can be considered as a probability
measure on $(\E, \F_{\E})$ as remarked earlier. 
All above functions are measurable in their arguments with respect to the
given $\sigma$-algebras, which is an easy application of the measurability
parts of Fubini's theorem and Campbell's theorem 
(see \cite{MKM}, Proposition 5.1.2. for example). 
Hence both $\pi_n$ and $\alpha_n$ are probability kernels. 
By the above definitions and by \eqref{azerlsup} for every  
$D \in \F_{\Y}$ and $\bar{Y} \in \Y$  one has the decomposition 
\begin{equation} \label{zersup}
\begin{split}
&\gmu_{\Lambda_n}(D \cap \X_n|\bar{Y})\\  
&= \;  \frac{1}{Z^U_{\Lambda_n}(\bar{Y})} 
 \int \nu^o(dX) \int_{S^{X_{\Lambda_n}}} d\sigma_{X_{\Lambda_n}} \;
 e^{ -H^U_{\Lambda_n}((X_{\Lambda_n},\sigma)\ync)} 
 \; 1_{(D\cap \X_n)}((X_{\Lambda_n},\sigma)\ync )\\
&= \; \int \gmu_{\Lambda_n}^o(dX|\bar{Y}) \; 1_{\X_n}(X_{\Lambda_n}\ync^o) \; 
 \int \pi_n(dA|X_{\Lambda_n},\ync) \; \alpha_n(D|A,X_{\Lambda_n},\ync),
\end{split}
\end{equation}
where $\gmu_{\Lambda_n}^o(.|\bar{Y}) 
:= \gmu_{\Lambda_n}(.|\bar{Y}) \circ o^{-1}$.   
Now we want to examine the percolation process given by 
$\pi_n$. So let $n \in \N$, $Y \in \Y$ and $E_n := E(Y^o,n)$. 
$\pi_n(.|Y^o_{\Lambda_n},Y_{\Lambda_n^c})$ has its whole weight 
on the countable set of finite subsets  $A \subset E_n$, 
but this measure shows a strong dependence of different bonds. 
Fortunately, this measure is stochastically dominated ($\preceq$)  
by a Bernoulli measure, where the order on
the underlying space $\E_{E_n}$ is given by the inclusion. 
This stochastic domination will be an
important tool for evaluating bond probabilities. For a definition of
stochastic domination see \cite{GHM}, for example.\\  
More precisely, for given $X \in \X$ let $\pi(.|X)$ be the
Bernoulli measure on  $(\E_{E(X)},\Borel_{E(X)})$ with bond
probabilities  $\epsilon_{x_1x_2} :=  J(x_1-x_2) \epsilon$ 
for $x_1x_2 \in E(X)$. Note that  $0 \le \epsilon_{x_1x_2} \le 1$ for
all bonds $x_1x_2 \in E(X)$, which is a consequence of the condition on  
$\epsilon$ in \eqref{bedepssup}, and even $0 < \epsilon_{x_1x_2}$ for 
all $x_1x_2 \in E(X,n)$. Again $\pi(.|X)$ 
can be considered as a probability measure on $(\E, \F_{\E})$, and
indeed is a probability kernel. We now have
\begin{lemma} \label{ledominierungsup}
For all $n \in \N$ and $Y \in \Y$,  
\begin{equation} \label{dominierungsup}
\pi_n(.|Y^o_{\Lambda_n},Y_{\Lambda_n^c}) \; \preceq \; \pi(.|Y^o).
\end{equation}
\end{lemma}

%%%%%%%%%%%%%%%%%%%%%%%%%%%%%%%%%%%%%%%%%%%%%%%%%%%%%%%%%%%%%%%%%%%%

\subsection{Deforming the spin transformations}

For a configuration of positions $X \in \X$ and a bond set 
$A \subset E(X)$ let $\; \stackrel{A}{\longleftrightarrow} \; := \;  
\stackrel{A,X}{\longleftrightarrow} \; $ be the equivalence relation 
on $X$ such that for all $x_1, x_2 \in X$ we have 
$x_1 \stackrel{A,X}{\longleftrightarrow} x_2$ iff either $x_1 = x_2$
or there is a finite path in $X$ joining $x_1$ and $x_2$ and using
bonds in $A$ only. 
For $x_1 \ne x_2  \in X$, the inequality 
\begin{equation} \label{Abpsup}
 \pi(x_1 \stackrel{.}{\longleftrightarrow} x_2 | X) \quad 
\le \quad   \sum_{m \ge 1} \quad \sideset{}{^{\neq}}
 \sum_{\substack{x'_0,...,x'_m \in X: \\
 x'_0=x_1,x'_m=x_2}} \epsilon^m \; \prod_{i=1}^{m}  
 J(x'_i - x'_{i-1})
\end{equation}
is an easy consequence of the above definition.
For a configuration $X \in \X$, a bond set $A \subset E(X)$ and a
point $x \in X$ let 
\begin{displaymath}
C_{A,X}(x) :=\{ x' \in  X: x \stackrel{A}{\longleftrightarrow} x' \}   
\end{displaymath}
be the percolation cluster of $x$ in $(X,A)$.    
Furthermore we want to consider the range of clusters, so for 
$x \in X$ and $\Lambda \in \Borel^2_b$ let 
\begin{displaymath}\begin{split}
r_{A,X}(x) \quad &:= \quad \sup\{ \|x'\| : x' \in C_{A,X}(x) \} 
\quad \text{ and } \\
r_{A,X}(\Lambda) \quad &:= \quad 
\left \{
\begin{aligned}
&\max\{r_{A,X}(x'): x'\in \Lambda \cap X \}  &&\text{ for } \Lambda
\cap X \ne \emptyset \\
&0 &&\text{ for } \Lambda \cap X = \emptyset.
\end{aligned} \right.\\
\end{split}
\end{displaymath}
Obviously $\|x\| \le r_{A,X}(x) \le \infty$ and
$r_{A,X}(\Lambda) \le \infty$. Now we have an estimate for the range
of the cluster of the given set $\Lambda_{n'}$, where $n'$ is the natural
number fixed in section \ref{constantssup}. 
\begin{lemma} \label{lereichsup}
There exists an integer $R > n'$ and a set $\Xrn \in \F_{\X}$ such that  
$\mu(\Xrn) \ge 1-2\delta$ and, for every $Y \in \Xrn$ and $n \ge n'$,  
\begin{equation} \label{reichsup}
\pi_n(\{ A: r_{A,Y^o}(\Lambda_{n'}) \ge R \} \; 
  |\; Y_{\Lambda_n}^o, Y_{\Lambda_n^c}) \;
\le \; \delta. 
\end{equation}
\end{lemma}
From now on let an integer $R \ge 2$ with the above property be
fixed.
In order to construct the spin deformation we define the functions 
$q: \R \to \R$, $Q: \R \to \R$, $r: \R \times \R_+ \to \R$ and 
$\tau_{n}: \R^2 \to S^1$ for $n > R$  by  
\begin{displaymath}
\begin{split}
q(s) \; &:= \; 1_{\{s \le 2\}} \; + \; \frac{1}{s \log(s)} 1_{\{s > 2\}},
\hspace{ 2.5 cm}Q(k) := \int_0^k q(s) ds, \\
r(s,k) \; &:= \; 1_{\{s \le 0\}} \; 
 + \; \int_s^k \frac{q(s')}{Q(k)} ds' 1_{\{0 < s < k\}}, \qquad   
\tau_{n}(x) := \tau \cdot r(\|x\|-R,n-R).
\end{split}
\end{displaymath}
\begin{lemma} \label{leeigQrtausup}
For all $n>R$ and $x,x' \in \R^2$ such that  
$\|x'\| \ge \|x\|$ we have  
\begin{equation} 
\label{taugegenqsup}
0 \quad \le \quad \tau_n(x) - \tau_n(x') \quad
\le \quad \tau \|x-x'\| \frac{q(\|x\| -  R)}{Q(n-R)},
\end{equation}
$\lim_{n \to \infty}Q(n) = \infty$,   
\begin{equation} \label{tauinnenaussensup}
\tau_{n}(x) = \tau \; \text{ for } \; \|x\| \le R \qquad 
\text{ and } \qquad \tau_{n}(x) = 0 \; \text{ for } \; \|x\| \ge n.
\end{equation}
\end{lemma}
However, what we really need here is a spin deformation which is constant on 
points joined by a bond of a given set $A$. So, for $n \in \N$, 
$X \in \X$ and  $A \subset E(X,n)$ we define 
$\tau^{X,A}_n: X \to S^1$ by
\begin{displaymath}
\tau^{X,A}_n(x) \; :=  \; 
\min \, \{ \,\tau_{n}(x') : x'\in X \text{ and } x
\stackrel{A}{\longleftrightarrow} x' \, \}. 
\end{displaymath}
This spin deformation can be seen to be measurable in $x,X$ and $A$ with respect
to the given $\sigma$-algebras using Campbell's theorem. 
Because of \eqref{tauinnenaussensup} we have $\tau_{n}(x') = 0$ for 
$\|x'\| \ge n$, so the minimum is attained at some point 
$t_A(x) \in X$ ($t_A(x) := x$ for $\|x\| \ge n$). By construction we have
\begin{equation} \label{gleichfallskongruentsup}
\begin{split}
&\|t_A(x)\| \ge \|x\|, \quad
t_A(x) \stackrel{A}{\longleftrightarrow} x, \quad 
\tau^{X,A}_n(x) = \tau_n(t_A(x)) \quad \alle x \in X\\
&\text{ and } \qquad \tau^{X,A}_n(x) = \tau^{X,A}_n(x') \quad 
\alle x,x' \in X
\text{ such that } x \stackrel{A}{\longleftrightarrow} x'. 
\end{split}
\end{equation}

%%%%%%%%%%%%%%%%%%%%%%%%%%%%%%%%%%%%%%%%%%%%%%%%%%%%%%%%%%%%%%%%%%%%%%

\subsection{Proof of Theorem \ref{symsup}}

In order to simplify notation, for $n \in \N$, $X \in \X$ and 
$E_n := E(X,n)$ let
$f_{n,X}: \E_{E_n} \to \R$ be defined by   
\begin{equation}
f_{n,X}(A) \; := \; \sum_{xx' \in E_n}
J(x-x')(\tau^{X,A}_n(x) -\tau^{X,A}_n(x'))^2.
\end{equation}
\begin{lemma} \label{lepirestsup}
There exists an integer $n > R$ and a set of configurations
$\Xbrn \in \F_{\X}$ such that   
$\mu(\Xbrn) \ge 1-\delta$ and, for every $Y \in \Xbrn$,   
\begin{equation} \label{pirestsup}
\pi_n \Big( f_{n,Y^o} \ge \frac{2}{\|\bar{V}''\|} 
 \; \Big| \;   Y_{\Lambda_n}^o,Y_{\Lambda_n^c} \Big) 
\; \le \; \delta.
\end{equation}
\end{lemma}
Let such an $n$ be fixed for the rest of the proof, let 
$\Xd := \Xbrn \cap \Xrn \cap \Xn$ be the set of good configurations of
positions, and for $X \in \X$ let  
\begin{displaymath}
\alg_{n,X} \; := \; \Big \{ A \subset E(X,n): r_{A,X}(\Lambda_{n'}) < R,\; 
f_{n,X}(A) < \frac{2}{ \|\bar{V}'' \| } \Big \}
\end{displaymath}
be the set of good bond sets. 
\begin{lemma} \label{letaylorsup}
For every $Y \in \Xd$ and $A \in \alg_{n,Y^o}$ we have 
\begin{align} 
\label{deltasup}
&\mu(\Xd) \; \ge \; 1 - 3\delta \quad \text{ and } \quad
\pi_n(\alg_{n,Y^o} \Bed Y_{\Lambda_n}^o,Y_{\Lambda_n^c}) \; 
\ge \; 1- 2\delta,\\
\label{tauAsup}
&\tau^{Y^o,A}_n(x) \; = \; \tau  \;\; \alle x \in Y^o_{\Lambda_{n'}}
\quad \text{ and } \quad 
\tau^{Y^o,A}_n(x) \; = \; 0 \;\; \alle x \in Y^o_{\Lambda_n^c}, \\ 
\label{taylorsup}
&\frac{e}{2} e^{ -H^{\bar{U}}_{\Lambda_n}((\tau^{Y^o,A}_n)^{-1}Y)} \; 
 + \; \frac{e}{2} e^{ -H^{\bar{U}}_{\Lambda_n}(\tau^{Y^o,A}_nY)} \quad \ge
 \quad  e^{ -H^{\bar{U}}_{\Lambda_n}(Y)}.
\end{align}
\end{lemma}
All these facts together imply  
\begin{lemma} \label{finishsup}
For the integer $n$ and the set $\Xd$ we have 
\begin{equation} \label{georgiisup}
\frac{e}{2} \gmu_{\Lambda_n}(\tau^{-1}B \cap \Xd |\bar{Y}) \; + \; 
\frac{e}{2} \gmu_{\Lambda_n}(\tau B \cap \Xd |\bar{Y}) \quad 
\ge \quad \gmu_{\Lambda_n}(B \cap \Xd |\bar{Y}) \; - \; 2\delta.
\end{equation}
\end{lemma}
Now  integrating \eqref{georgiisup} - using property
\eqref{gibbsaequsup} of $\mu$ and \eqref{deltasup} -  yields 
\begin{equation*}
\frac{e}{2} \mu(\tau^{-1}B) \; + \; 
\frac{e}{2}\mu(\tau B) \quad \ge \quad \mu(B) \; - \; 5 \delta
\end{equation*}
for arbitrary $\mu \in \G(U)$, $\tau \in G$, $n' \in \N$,$B \in \F_{\Y,
  \Lambda_{n'}}$ and $\delta > 0$. Letting $\delta \to 0$
the assertion of the theorem follows by using results from the
general theory of Gibbs measures, see \cite{G}, chapter 9.1,
proposition (9.1) for example.

%%%%%%%%%%%%%%%%%%%%%%%%%%%%%%%%%%%%%%%%%%%%%%%%%%%%%%%%%%%%%%%%%%%%%%
%%%%%%%%%%%%%%%%%%%%%%%%%%%%%%%%%%%%%%%%%%%%%%%%%%%%%%%%%%%%%%%%%%%%%%

\section{Proofs of the lemmas}

\subsection{Property of the correlation function: 
Lemma \ref{lekorabsup}}

Let $U$ be a superstable and lower regular potential, 
$\mu \in \G(U)$ a tempered Gibbs measure, $\xi \in \R$ a
correlation bound, $m \ge 0$ an integer and $f: (\R^2)^m \to \R_+$ 
a measurable function. 
The Poisson point process $\nu$  satisfies for every 
measurable $g: \Y_{\Lambda_n} \to \R_+$ 
\begin{displaymath}
\begin{split}
\int &\nu(dY)  \sideset{}{^{\neq}}\sum_{x_1,...,x_m \in 
    Y^o_{\Lambda_N}} f(x_1,...,x_m) \; g(Y) \\ 
&= \; z^m \int_{\Lambda_N^m}dx_1 ... dx_m 
  \int_{E^m} d\sigma_1... d\sigma_m  \; f(x_1,...,x_m) \;
  \int \nu(dY')\; g((X,\sigma)_m Y'),
\end{split}
\end{displaymath}
where $ (X,\sigma)_m := \{ (x_i,\sigma_i): 1 \le i \le m \}$. 
Using this equality, the characterisation of Gibbs measures 
\eqref{gibbsaequsup}, the definition of the conditional Gibbs 
distribution and the definition of the correlation function we get 
\begin{displaymath}
\begin{split}
&\int \mu(dY) \sideset{}{^{\neq}}\sum_{x_1,...x_m \in Y^o_{\Lambda_N}} 
 f(x_1,...,x_m)\\
&= \;  \int \mu(d\bar{Y}) \; \frac{1}{Z^U_{\Lambda_{N}}(\bar{Y})}
   \; \int \nu (dY)  
   \sideset{}{^{\neq}} \sum_{x_1,...x_m \in Y^o_{\Lambda_N}}  
   f(x_1,...,x_m) \; 
    e^{- H^U_{\Lambda_{N}}(Y_{\Lambda_N}\bar{Y}_{\Lambda_{N}^c})}\\   
&= \; \int_{\Lambda_N^m}dx_1 ... dx_m \int  d\sigma_1... d\sigma_m
   \; f(x_1,...,x_m) \; z^m \; \rho^{U,\mu} ((X,\sigma)_m)\\ 
&\le \; (z \xi)^m \; \int_{\Lambda_N^m} dx_1 ... dx_m  
  \; f(x_1,...,x_m),  
\end{split} 
\end{displaymath} 
where we have used the bound \eqref{korbeschrsup} on the correlation function 
in the last step. Letting $N \to \infty$ the assertion \eqref{korabsup}
follows from the monotone limit theorem.

%%%%%%%%%%%%%%%%%%%%%%%%%%%%%%%%%%%%%%%%%%%%%%%%%%%%%%%%%%%%%%

\subsection{Convergence of energy sums: Lemma \ref{umordsup}}

Let $n \in \N$. For every $X \in \X$ we have 
\begin{displaymath}
\sum_{x_1x_2 \subset E(X,n)} J(x_1-x_2) \quad 
\le \quad 
\sideset{}{^{\neq}}\sum_{x_1, x_2 \in X} 
1_{\{x_1 \in \Lambda_n\}} \; J(x_1-x_2), \quad \text{ so } 
\end{displaymath}
\begin{displaymath}
\int \mu(dY) \sum_{x_1x_2 \subset E(Y^o,n)} J(x_1-x_2) \quad
\le \quad (z \xi)^2 \int dx_1 dx_2  
 \; 1_{\{x_1 \in \Lambda_n\}} \;J(x_1-x_2) 
\end{displaymath}
by lemma \ref{lekorabsup}, and the right hand side of the last 
inequality is at most $\;c_J (2nz\xi)^2\; < \infty$.
So the assertion is true for  
\begin{displaymath}
\Xn \; := \; \Big \{ \;  X \in \X: \;
\sum_{x_1x_2 \subset E(X,n)} J(x_1-x_2)  \; < \; \infty \; \Big \}. 
\end{displaymath}

%%%%%%%%%%%%%%%%%%%%%%%%%%%%%%%%%%%%%%%%

\subsection{Stochastic domination:  Lemma \ref{ledominierungsup}}

A general sufficient condition for stochastic domination in a situation
like the one considered is given by R. Holley (see \cite{H}
e. g.). The result is the following:

\begin{lemma} \label{ledomikritsup}
Let $Z = \{ e_1, e_2, ...\}$ be a countable set, 
$(\epsilon_e)_{e \in Z}$ a familiy of reals in $]0,1]$,  
$\Borel_{Z}$ the Bernoulli-$\sigma$-algebra on $\p(Z)$, and let $\A$
and $\A_{\epsilon}$ be random variables with values in 
$(\p(Z), \Borel_{Z})$ such that $\A_{\epsilon}$ is a
Bernoulli process with bond probabilities $\epsilon_e$, 
and for every $e \in Z$ we have $P( e \in \A | \A \weg e) \le \epsilon_e$
a. s. . 
Then $\V(\A) \preceq \V(\A_{\epsilon})$. 
\end{lemma} 

\Bew
Let all assumptions of the lemma hold. First we consider the finite 
sets $Z^{(n)} := \{e_1, ... , e_n\}$ and let  $\A^{(n)},
\A^{(n)}_{\epsilon}$ be the restrictions of $\A, \A_{\epsilon}$ to 
$Z^{(n)}$, i. e. $\A^{(n)} =  \A \cap Z^{(n)}$ and 
$\A^{(n)}_{\epsilon} =  \A_{\epsilon} \cap Z^{(n)}$. 
For any $n \in \N$ and $e \in Z^{(n)}$ we have
$P( e \in \A^{(n)} | \A^{(n)} \weg e) \le \epsilon_e$ a. s. , 
which is a straightforward consequence of 
$P( e \in \A | \A \weg e) \le \epsilon_e$ a. s.  and the
properties of conditional probabilities. 
Now the criterion of R. Holley (as presented in \cite{GHM}, Theorem
4.8., for example) gives  
$\V(\A^{(n)}) \preceq \V(\A^{(n)}_{\epsilon})$. If
$\V(\A^{(n)})$ and $\V(\A^{(n)}_{\epsilon})$ are considered as measures on
$(\p(Z),\Borel_{Z})$ we observe that
\begin{equation*} 
\V(\A^{(n)}) \to \V(\A) \quad \text{ and } \quad
\V(\A^{(n)}_{\epsilon}) \to 
\V(\A_{\epsilon}) \quad \text{ weakly } \quad \text{ as } n \to
\infty.  
\end{equation*}
As stochastic domination is preserved under weak limits 
(see \cite{GHM}, Cor. 4.7.,  for example) we get  
$\V(\A) \preceq \V(\A_{\epsilon})$. \qed\\

Now, turning to the proof of lemma \ref{ledominierungsup}
let $n \in \N$, $Y \in \Y$ and $E_n := E(Y^o,n)$. In order to show 
that $\pi_n(.|Y^o_{\Lambda_n},Y_{\Lambda_n^c}) \; \preceq \;
\pi(.|Y^o)$  we may consider both measures as measures on
$(\E_{E_n},\Borel_{E_n})$. We also may assume that 
$Y^o \in \Xn$ and $\W_n(Y_{\Lambda_n},Y^o_{\Lambda_n^c}) > 0$. 
By lemma \ref{ledomikritsup} it is sufficient to show that, 
for every bond $x_1x_2 \in E_n $ and every finite bond set 
$D \subset E_n \weg \{x_1x_2\}$, 
\begin{displaymath}
\pi_n(\{x_1x_2\} \cup D \Bed Y^o_{\Lambda_n},Y_{\Lambda_n^c}) \; 
 \le \; \epsilon_{x_1x_2} \; \pi_n(\{D,\{x_1,x_2\} \cup D\} 
\Bed Y^o_{\Lambda_n},Y_{\Lambda_n^c}) .
\end{displaymath}
(Here we have used that the whole weight of
$\pi_n(.|Y^o_{\Lambda_n},Y_{\Lambda_n^c})$ is on the 
countable set of finite bond sets.)  
So let $x_1x_2 \in E_n$ and $D \subset E_n \weg \{x_1x_2\}$ be finite. 
By the definition of $\pi_n$ the last inequality is equivalent to 
\begin{displaymath}
\begin{split} 
 &\int d\sigma_{Y^o_{\Lambda_n}} \; 
\w_n(D,(Y^o_{\Lambda_n},\sigma)Y_{\Lambda_n^c}) \; 
\Big[ \; \epsilon_{x_1x_2}\\
&\qquad  + 
 (\epsilon_{x_1x_2} - 1) \Big(e^{ J(x_1-x_2)
v(\sigma_{x_1}((Y^o_{\Lambda_n},\sigma)\ync)
-\sigma_{x_2}((Y^o_{\Lambda_n},\sigma)\ync) )} -1 \Big)\; \Big] \;\ge\; 0.
\end{split}
\end{displaymath}
But since $0 < \epsilon_{x_1x_2} \le 1$ and $0 < v < \epsilon$, 
the term in the brackets is at least  
\begin{displaymath}
\epsilon_{x_1x_2} +
(\epsilon_{x_1x_2}-1)(e^{\epsilon_{x_1x_2}} - 1) \quad \ge  \quad 0,
\end{displaymath}
which completes the proof of the Lemma \ref{ledominierungsup}.

%%%%%%%%%%%%%%%%%%%%%%%%%%%%%%%%%%%%%%%%%%%%%%%%%%%%%%

\subsection{Cluster bounds: Lemma \ref{lereichsup}}

Let $n \ge n'$ be a fixed integer. For a given configuration 
$X \in \X$ and a bond set $A \in E(X,n)$ we consider the cardinality
of the cluster of points from  $\Lambda := \Lambda_{n'}$, 
which is defined by
 \begin{displaymath}
C_{\Lambda}(A) \quad := 
\quad |\bigcup_{x \in X_{\Lambda}}  C_{A,X}(x)|. 
\end{displaymath}
For all $X \in \X $ we have the estimate 
\begin{displaymath}
\begin{split}
\int& \pi(dA|X)\; C_{\Lambda}(A) \quad 
\le \quad  \int \pi (dA|X) \; 
 \sum_{ x \in X_{\Lambda}} \sum_{x' \in X} 
1_{\{x \stackrel{A}{\longleftrightarrow} x'\}} \\
&=  \sum_{ x \in X_{\Lambda}}  \sum_{x' \in X}
 \pi(x \stackrel{.}{\longleftrightarrow} x'|X)\\  
&\le\;  \sum_{m \ge 0} \epsilon^m 
 \sideset{}{^{\neq}}\sum_{x_0,...,x_m \in X} 
 1_{x_0 \in \Lambda} \;  \prod_{i=1}^{m} J(x_i - x_{i-1}) \quad =: \quad 
f(X),
\end{split}
\end{displaymath}
where we have used \eqref{Abpsup}. By Lemma \ref{lekorabsup} we have  
\begin{equation*}
\begin{split}
\int \mu(dY) &f(Y^o) \quad 
\le \; \sum_{m \ge 0} \; \epsilon^m \; (z\xi)^{m+1} \; 
   \int dx_0 ... dx_{m} \; 1_{x_0 \in \Lambda} \;  
   \prod_{i=1}^{m} J(x_{i} - x_{i-1})\\ 
&\le \; z \xi(2n')^2 \;  \sum_{m \ge 0} \; (z \xi \epsilon c_J)^m 
\; =: \; c \; < \; \infty
\end{split}
\end{equation*}
due to \eqref{bedepssup}. Letting  
\begin{equation*}
\Xrn' \quad := \quad \big \{\;  X \in \X: \; f(X)  
\; \le \;  \frac{c}{\delta} \; \big\} 
\end{equation*}
we get $\mu(\Xrn') \ge 1 - \delta$ from Chebyshev's inequality, 
and for any $X \in \Xrn'$ we have again by Chebyshev's inequality that 
\begin{displaymath}
\pi \Big(C_{\Lambda} > \frac{2c}{\delta^3}\; \Big| \; X \Big) 
 \quad \le \quad \frac{\delta^3}{2c}  
 \int \pi(dA|X)\; 
 C_{\Lambda}(A) \quad \le \quad \frac{\delta^2}{2}.
\end{displaymath}
Now let $n \ge n'$, $R > n'$ and $X \in \Xrn'$. Then, by
the above estimate,  
\begin{displaymath}
\begin{split}
\pi&(r_{.,X}(\Lambda) \ge R \Bed X)\\
&\le \; \pi \Big( C_{\Lambda} > \frac{2c}{\delta^3} \; \Big| \;  X \Big) 
 \; + \;  \pi \Big( C_{\Lambda} \le \frac{2c}{\delta^3} \; , \; 
r_{\;.\;,X}(\Lambda) \ge R \; \Big| \; X \Big)\\
&\le \; \frac{\delta^2}{2} \; + \;
 \pi \Big( \big \{ A:\gibt 1 \le m \le \frac{2c}{\delta^3} \; 
 \gibt \text{ distinct }x_0,..., x_m \in X: \\
& \quad  \gibt 1 \le j \le m:
x_0 \in \Lambda,\; \|x_j - x_{j-1}\| \ge
\frac{\scriptstyle (R-n')\delta^3}{ \scriptstyle 2c}, \; 
x_{i-1}x_i \in A \alle i \big\} \; \Big| \; X\Big)\\ 
&= \; \frac{\delta^2}{2} \; + \; \sum_{m \ge 1} \; \sum_{j=1}^m
 \; \; \sideset{}{^{\neq}} \sum_{x_0,..., x_m \in X}
 1_{\{x_0 \in \Lambda , \|x_j - x_{j-1}\| \ge \frac{(R-n')\delta^3}{2c}\}}
 \; \epsilon^m  \prod_{i=1}^m J(x_i-x_{i-1})\\
&=: \; \frac{\delta^2}{2} \; + \;  f_{R}(X),
\end{split}
\end{displaymath}
and Lemma \ref{lekorabsup} yields 
\begin{displaymath}
\begin{split}
\int \mu(dY) f_{R} (Y^o) \; 
&\le \; \sum_{m \ge 1} \; \epsilon^m \; \sum_{j=1}^m \; (z\xi)^{m+1} 
  \;  \int dx_0 ... dx_{m} \; \Big[ \\
&\hspace{2cm}  1_{\{x_0 \in \Lambda , 
       \|x_{j} - x_{j-1}\| \ge \frac{(R-n')\delta^3}{2c}\}}\; 
 \prod_{i=1}^m J(x_{i}-x_{i-1}) \; \Big]\\
&\le \; z\xi \;  \sum_{m \ge 1} \; (z\xi\epsilon)^m \; m \; 
  (2n)'^2 \; c_J^{m-1} \; 
  c \Big( \frac{\scriptstyle (R-n')\delta^3}{\scriptstyle 2c} \Big). 
\end{split}
\end{displaymath}
In the last step, the integrals have been estimated backwards from
$x_{m}$ to  $x_0$, where integration over $x_j$ gives the constant 
$c\big(\frac{(R-n')\delta^3}{2c}\big)$ defined in \eqref{crsup}. As 
$\lim_{R \to \infty} c \big(\frac{(R-n')\delta^3}{2c}\big ) = 0$ 
and the sum over $m$ is finite by condition \eqref{bedepssup}, we can fix
an $R > n'$ such that   
\begin{equation} \label{intintcebsup}
\int \mu(dY) \Big( \; \frac{\delta^2}{2} \; +
 f_{R} (Y^o) \; \Big) \quad \le \quad \delta^2.
\end{equation}
Now let  
\begin{equation*}
\Xrn'' \; := \; \Big \{ \;  X \in \X:
\; \frac{\delta^2}{2} \; + \; f_{R} (X) \;  \le \; 
\delta \; \Big \}  
\end{equation*}
and $\Xrn \; := \; \Xrn'' \cap \Xrn'$, then by Chebyshev's inequality 
and \eqref{intintcebsup} we have 
$\mu(\Xrn'') \ge 1 - \delta$, and hence $\mu(\Xrn) \ge 1 - 2\delta$. 
For every $Y  \in \Xrn$ the event  
$\{ A: r_{A,Y^o}(\Lambda) \ge R \}$ is increasing, so by stochastic
domination \eqref{dominierungsup} we have  
\begin{displaymath}
\begin{split}
\pi_n&(\{ A: r_{A,Y^o}(\Lambda) \ge R \} \; |\; Y_{\Lambda_n}^o,
 Y_{\Lambda_n^c})
\quad \le \quad \pi_{n,\epsilon}
 (\{ A: r_{A,Y^o}(\Lambda) \ge R \} \; |\; Y^o) \\
&\le \quad \frac{\delta^2}{2} \; + \; f_{R}(Y^o) \quad \le \quad
\delta. 
\end{split}
\end{displaymath}

%%%%%%%%%%%%%%%%%%%%%%%%%%%%%%%%%%%%%%%%%%%%%%%%%%%%%%%%%%%%%%%%%%%%

\subsection{Properties of $\tau_n$ and $Q$: Lemma \ref{leeigQrtausup}}

\eqref{tauinnenaussensup} is evident from the
definition of $\tau_{n}$, and 
$\lim_{n \to \infty} Q(n) = \infty$ is a consequence of $\log \log n
\le Q(n)$ for $n \ge 2$.  
For \eqref{taugegenqsup} let $x,x' \in \R^2$ such that 
$\|x'\| \ge \|x\|$. The left inequality is trivial and for the right
inequality we may assume that $\|x'\| > R$ and $\|x\| < n$ because of 
\eqref{tauinnenaussensup}. Hence 
\begin{displaymath}
\begin{split}
r&(\|x\|-R,n-R) - r(\|x'\|-R,n-R) \quad 
= \quad  \int_{\max\{R,\|x\|\}}^{\min\{\|x'\|, n\}}
\frac{q(s'-R )}{Q(n-R)} ds'\\ 
&\le \quad  (\|x'\|-\|x\|) \; \frac{q(\|x\| - R)}{Q(n-R)} \quad 
 \le \quad  \|x'-x\| \; \frac{q(\|x\|-R)}{Q(n-R)},
\end{split}
\end{displaymath}
where we have used the monotonicity of $q$ and the triangle
inequality. Now \eqref{taugegenqsup} follows immediately.

%%%%%%%%%%%%%%%%%%%%%%%%%%%%%%%%%%%%%%%%%%%%%%%%%%%%%%%%%%%%%%%%%%%%%

\subsection{Probability of bad bond sets: Lemma
\ref{lepirestsup}} 

First of all we state two easy facts. First, 
\begin{equation}\label{unglquadsup}
\|x_m-x_0\|^2  \; 
\le \;  m \prod_{i=1}^m (\|x_i - x_{i-1}\|^2 +1) \quad 
\alle m \ge 1, x_0,...,x_m \in \R^2, 
\end{equation}
by the triangle inequality and the arithmetic-quadratic mean inequality. 
Secondly, 
\begin{equation}\label{intqsup}
\int_{\Lambda_n} dx  \; q(\|x\| -  R)^2 \quad 
\le \quad 8(R+3)^2 \;+ \; 8R Q(n-R) \qquad \alle n \ge R,    
\end{equation}
which is obtained by the substitution $t := \|x\|$: 
\begin{displaymath}
\begin{split}
\int_{\Lambda_n} &dx  \; q(\|x\| -  R)^2 \quad 
\le \quad \int_0^{R+3} dt \; 8t \; 
  + \; \int_3^{n-R} dt \; 8(t+R) q(t)^2 \\
&\le  \; 8(R+3)^2 \;+ \; 8R \int_0^{n-R} q(t)dt \quad
=  \quad 8(R+3)^2 \;+ \; 8R Q(n-R),     
\end{split}
\end{displaymath}
where we have used in the first step that 
$q(t) \le 1 \alle t \in \R$, and in the second step that 
$t+R \le tR$ for $t,R \ge 2$, and $tq(t) \le 1 \alle t \ge 3$.\\

Now for the proof of Lemma \ref{lepirestsup} let $n > R$ and 
$Y \in \Y$ be arbitrary. Using the arithmetic-quadratic mean 
inequality to estimate $(\tau^{X,A}_n(x)-\tau^{X,A}_n(x'))^2 $ we get 
\begin{displaymath}
\begin{split}
f_{n,Y^o}(A) \;
&\le \; 6 \sum_{x,x' \in Y^o} 1_{\{x \neq x'\}} \;  
 J(x-x') \; (\tau_n(t_A(x))-\tau_n(x))^2  \\ 
&\qquad + 3 \sum_{x,x' \in Y^o}
 1_{\{\|x\| \le \|x'\|\}}\; 
 J(x-x') \; (\tau_n(x)-\tau_n(x'))^2. 
\end{split}
\end{displaymath}
Substituting $z := t_A(x)$ and introducing $1_{\{z = t_A(x)\}}$ in
the first sum we need only consider $z \in Y^o$ such that 
$\|x\| \le \|z\|$ and $x \neq z$. By distinguishing the cases 
$z \neq x,x'$ and $z = x'$ and by using $\{A:t_A(x) = z \} \; 
\subset \; \{A: x \stackrel{A}{\longleftrightarrow} z \}$ we can estimate
the expectation value of $f_{n,Y^o}$ by 
\begin{displaymath}
\begin{split}
\int&\pi_{n}(dA | Y^o_{\Lambda_n},Y_{\Lambda_n^c}) \; f_{n,Y^o}(A)\\
&\le \; 6 \sideset{}{^{\neq}}\sum_{x,x',z \in Y^o}
  1_{\{\|x\| \le \|z\|\}} \; J(x-x') \; (\tau_n(z)-\tau_n(x))^2 \; 
 \pi_{n}( x \stackrel{ . }{\longleftrightarrow} z
  |Y^o_{\Lambda_n},Y_{\Lambda_n^c})\\  
&\quad +  9 \sideset{}{^{\neq}}\sum_{x,z \in Y^o}  
 1_{\{\|x\| \le \|z\|\}}  \; J(x-z) \; (\tau_n(x)-\tau_n(z))^2 
\end{split}
\end{displaymath}
Next we use the stochastic domination \eqref{dominierungsup} for the
increasing events ${x \stackrel{ . }{\longleftrightarrow} z}$ 
to estimate 
$\pi_{n}( x \stackrel{ . }{\longleftrightarrow} z
  |Y^o_{\Lambda_n},Y_{\Lambda_n^c})$, and we use \eqref{taugegenqsup} from
  Lemma \ref{leeigQrtausup}, noting that  $\tau_n(x) = 0 = \tau_n(z)$ for 
$n <\|x\| \le \|z\|$. So we get
\begin{displaymath}
\begin{split}
\int\pi_{n}&(dA | Y^o_{\Lambda_n},Y_{\Lambda_n^c}) \; f_{n,Y^o}(A)\\
&\le \; 6 \sideset{}{^{\neq}}\sum_{x,x',z \in Y^o} 
 1_{\{x \in \Lambda_n\}} \; J(x-x') \; \tau^2 \; \|x-z\|^2 \;  
 \frac{q(\|x\| -  R)^2}{Q(n-R)^2} \; 
 \pi ( x \stackrel{ . }{\longleftrightarrow} z | Y^o)\\ 
&\qquad + \;  9 \sideset{}{^{\neq}}\sum_{x,z \in Y^o}  
  1_{\{x \in \Lambda_n\}} \; J(x-z) \; \tau^2 \|x-z\|^2 
 \frac{q(\|x\| -  R)^2}{Q(n-R)^2} \\ 
&=: \quad \Sigma_1 (Y^o,n) \quad  
 + \quad  \Sigma_2 (Y^o,n).
\end{split}
\end{displaymath}
In order to deal with $\Sigma_1(Y^o,n)$ we distinguish the paths 
$x_0,...,x_m$ from $x$ to $z$ analogously to \eqref{Abpsup} and
distinguish the cases $x_j = x'$ and $x_j \neq x' \alle j$. Hence 
\begin{displaymath}
\begin{split}
\Sigma_1(Y^o,n) \quad
&\le \quad  6 \; \sum_{m \ge 1} \epsilon^m  \; 
 \sideset{}{^{\neq}} \sum_{x',x_0,...,x_m \in Y^o}  
 1_{\{x_0 \in \Lambda_n\}} \; J(x_0-x')\\ 
&\qquad \qquad  \tau^2 \; \|x_0-x_m\|^2 \;  
  \frac{q(\|x_0\| -  R)^2}{Q(n-R)^2} 
  \; \prod_{i=1}^{m} J(x_i - x_{i-1})\\
&\quad + \; 6 \; \sum_{m \ge 1} \epsilon^m \sum_{j=1}^{m-1} \quad 
 \sideset{}{^{\neq}}\sum_{x_0,...,x_m \in Y^o}  
 1_{\{x_0 \in \Lambda_n\}} \; J(x_0-x_j)\\ 
&\qquad \qquad  \tau^2 \; \|x_0-x_m\|^2 \;  
 \frac{q(\|x_0\| -  R)^2}{Q(n-R)^2} \; 
 \prod_{i=1}^{m} J(x_i - x_{i-1}).
\end{split}
\end{displaymath}
Applying Lemma \ref{lekorabsup} we thus find 
\begin{displaymath}
\begin{split}
\int &\mu(dY) \; \Sigma_1 (Y^o,n)\\
&\le \; 6 \; \sum_{m \ge 1} \; \epsilon^m \; (z\xi)^{m+1} 
   \int dx_0 ... dx_{m} \; \Big[ \;  1_{\{x_0 \in \Lambda_n\}} \;
   \tau^2 \; \|x_0-x_{m}\|^2 \; \frac{q(\|x_0\| -  R)^2}{Q(n-R)^2}\\ 
&\hspace{2cm} \prod_{i=1}^{m} J(x_{i} - x_{i-1}) 
  \Big( \; z\xi \int dx' J(x_0-x') \; + \; \sum_{j=1}^{m-1} 
 J(x_0 - x_j) \Big) \Big].  
\end{split}
\end{displaymath}
After applying \eqref{unglquadsup} to $\|x_0-x_{m}\|^2$ and
estimating the parentheses $(.)$ by $z \xi c_J m$ 
we evaluate the integrals backwards from $x_{m}$ to  $x_1$:
\begin{displaymath}
\begin{split}
\int &\mu(dY) \Sigma_1 (Y^o,n)\\
&\le \; 6 \sum_{m \ge 1} m^2 \epsilon^m \; c_J (z\xi)^{m+2} \;  
   (2c_J)^m \tau^2 \int dx 
 \; 1_{\{x \in \Lambda_n\}} \frac{q(\|x\| -  R)^2}{Q(n-R)^2}\\
&\le \; 6 \;  c_J (\xi z \tau)^2   
  \;[ \sum_{m \ge 1} m^2 (2\epsilon c_J \xi z)^m ] \;
 \frac{8(R+3)^2 \;+ \; 8R Q(n-R)}{Q(n-R)^2},   
\end{split}
\end{displaymath}
where we have used \eqref{intqsup} in the last step. The expectation
value of  $\Sigma_2$ can be treated similarly, and the estimates  
together give 
\begin{displaymath}
\begin{split}
\int \mu(&dY) \; \big[ \; \Sigma_1(Y^o,n) + \Sigma_2(Y^o,n)\; \big] \\
&\le \; \Big[ \; 6  \; \sum_{m \ge 1} m^2(2\epsilon c_J \xi z)^m \; + 
  \; 9  \Big] \;c_J \; (\tau z \xi) ^2 \;    
 \frac{8(R+3)^2 \;+ \; 8R Q(n-R)}{Q(n-R)^2}.   
\end{split}
\end{displaymath}
The sum over $m$ is finite by the choice of $\epsilon$
\eqref{bedepssup}, and because of $\lim_{k \to \infty} Q(k) = \infty$ 
the fraction on the right hand side can be made arbitrarily small  
choosing $n$ large enough. So there is an integer $n > R$ such that 
\begin{equation} \label{intintsup}
\int \mu(dY) \; \big[ \; \Sigma_1(Y^o,n) + \Sigma_2(Y^o,n)\; \big]
 \quad \le \quad \frac{2\delta^2}{\|V''\|}.
\end{equation}
Let
\begin{displaymath}
\Xbrn \; := \;  \Big\{ \; X \in \X: \; 
 \Sigma_1(X,n) + \Sigma_2(X,n) \; \le \; \frac{2\delta}{\|V''\|} \; \Big\}, 
\end{displaymath}
then we have found $n$ and $\Xbrn$ as desired, 
as by \eqref{intintsup} and Chebyshev's inequality
 we have $\mu(\Xbrn) \ge 1 - \delta$,  
and for every $Y \in \Xbrn$ we have by the definition of $\Xbrn$, $\Sigma_1$ and
$\Sigma_2$ and again by Chebyshev's inequality 
\begin{displaymath}
\pi_n\Big(\;f_{n,Y^o} \ge \frac{2}{ \|\bar{V}'' \|} 
     \; \Big| \; Y_{\Lambda_n}^o,Y_{\Lambda_n^c}\Big) \quad \le \quad \delta. 
\end{displaymath}

%%%%%%%%%%%%%%%%%%%%%%%%%%%%%%%%%%%%%%%%%%%%%%%%%%%%%%%%%%%%%%%%%%%%%%%%%

\subsection{Properties of good configurations and bond sets: 
Lemma \ref{letaylorsup}} 

Let $Y \in \Xd$, $A \in \alg_{n,Y^o}$ and $E_n := E(Y^o,n)$.
The inequalities \eqref{deltasup} then follow immediately from 
Lemma \ref{lereichsup} and Lemma \ref{lepirestsup}. 
\eqref{tauAsup} follows from \eqref{tauinnenaussensup} 
because $r_{A,Y^o}(\Lambda_{n'}) < R$.
For \eqref{taylorsup} we consider 
$\bar{V}$ as a $2\pi$-periodic function on $\R$.  
By the smoothness of $\bar{V}$ we can use a Taylor expansion 
to obtain for all $a,b \in \R$    
\begin{displaymath}
\bar{V}(a+b) + \bar{V}(a-b) 
- 2 \bar{V}(a) \; \le \; \| \bar{V}''\| b^2,
\end{displaymath}
where $\| \bar{V}''\| < \infty$, as $\bar{V}''$ is continuous on a
compact space. W.l.o.g. we may assume the right hand side of 
\eqref{taylorsup} to be positive, hence 
$| H^{\bar{U}}_{\Lambda_n}(Y) | < \infty$. So we have, introducing
$\eta_{x_1,x_2} := \sigma_{x_1}(Y) - \sigma_{x_2}(Y)$ and 
$\vartheta_{x_1,x_2} := \tau^{X,A}_n(x_1) - \tau^{X,A}_n(x_2)$, 
\begin{displaymath}
\begin{split}
&H^{\bar{U}}_{\Lambda_n}((\tau_n^{X,A})^{-1}Y)\;  
+ \;  H^{\bar{U}}_{\Lambda_n}(\tau_n^{X,A} Y) \;  
- \; 2H^{\bar{U}}_{\Lambda_n}(Y)\\ 
& \;= 
\sum_{x_1x_2 \in E_n} J(x_1-x_2) \; 
\left [\bar{V}(\eta_{x_1,x_2} - \vartheta_{x_1,x_2}) 
+ \bar{V}(\eta_{x_1,x_2} + \vartheta_{x_1,x_2}) 
 -   2  \bar{V}(\eta_{x_1,x_2}) \right]\\  
& \; \le \sum_{x_1x_2 \in E_n} J(x_1-x_2) \; \|\bar{V}''\| 
 \; \vartheta_{x_1,x_2}^2 \quad 
= \quad \|\bar{V}''\| \;  f_{n,Y^o}(A).
\end{split}
\end{displaymath}
By the convexity of the exponential function we conclude 
\begin{displaymath}
\begin{split}
\frac{e}{2}& e^{ -H^{\bar{U}}_{\Lambda_n}((\tau^{X,A}_n)^{-1}Y)}
 \; + \; \frac{e}{2} e^{ -H^{\bar{U}}_{\Lambda_n}(\tau^{X,A}_n Y)} \quad
\ge \quad  e^{1 - \frac{1}{2}
H^{\bar{U}}_{\Lambda_n}((\tau^{X,A}_n)^{-1}Y)
- \frac{1}{2}
H^{\bar{U}}_{\Lambda_n}(\tau^{X,A}_n Y)}\\
&\ge \; e^{1-\frac{\|\bar{V}''\|}{2} f_{n,Y^o}(A)}\cdot 
e^{-H^{\bar{U}}_{\Lambda_n}(Y)} \quad 
\ge \quad e^{- H^{\bar{U}}_{\Lambda_n}(Y)}. 
\end{split}
\end{displaymath}

%%%%%%%%%%%%%%%%%%%%%%%%%%%%%%%%%%%%%%%%%%%%%%%%%%%%%%%%%%%%%%%

\subsection{Inequality for the specifications: Lemma \ref{finishsup}}

By \eqref{zersup} it is sufficient to prove that for
every $Y \in \Xd$  we have   
\begin{displaymath}
\begin{split}
\int &\pi_n(dA|Y_{\Lambda_n}^o, Y_{\Lambda_n^c}) \; 
 \Big[ \; \frac{e}{2} \alpha_n
 (\tau^{-1}B | A,Y_{\Lambda_n}^o, Y_{\Lambda_n^c}) \\  
&\hspace{1cm}  + \; \frac{e}{2} 
 \alpha_n(\tau B | A,Y_{\Lambda_n}^o, Y_{\Lambda_n^c})  
 \; - \; \alpha_n(B|A,Y_{\Lambda_n}^o, Y_{\Lambda_n^c}) \Big] \; 
 + \; 2\delta  \quad \ge \quad 0, 
\end{split}
\end{displaymath}
and because of \eqref{deltasup} it suffices to show that, for every 
$Y \in \Xd$  and every finite $A \in \A_{n,Y^o}$ such that 
$\W_n(A, Y_{\Lambda_n}^o,Y_{\Lambda_n^c}) > 0$, 
the term in square brackets is nonnegative. 
By definition of $\alpha_n$, this will follow once we have shown that 
for every $X \in \X_{\Lambda_n}$, $Y \in \Y_{\Lambda_n^c}$  
and every finite $A \in \A_{n,XY^o}$ we have   
\begin{equation} \label{redsup}
\begin{split}
\int d\sigma'_X e^{ -H^{\bar{U}}_{\Lambda_n}(\yks)}\; 
 &\Big( \;  \frac{e}{2}1_{\tau^{-1}B} (\yks)
 + \frac{e}{2}1_{\tau B} (\yks)  - 1_{B} (\yks) \;  \Big)\\
& \cdot \prod_{x_1x_2 \in A}
 (e^{J(x_1-x_2)v(\sigma_{x_1}(\yks) - \sigma_{x_2}(\yks))}-1) 
 \quad \ge \quad 0,
\end{split}
\end{equation}
where we have used the notation $\yk := (X,\sigma)Y$. 
So let $X$, $Y$ and $A$ as above. The integral on
the left hand side of \eqref{redsup} can be split into the three 
parts $I_-,I_+$ and $I_0$  corresponding to the 
terms  $\tau^{-1}B$, $\tau B$ and $B$, and for any $x \in X$ 
we make the substitutions $\sigma_x := \sigma_x' + \tau^{X,A}_n(x)$ and 
$\sigma_x := \sigma_x' - \tau^{X,A}_n(x)$ in $I_-$ and $I_+$ respectively.
Because of \eqref{tauAsup} the spin transformation $\tau^{X,A}_n$ 
has no effect outside of $\Lambda_n$, so that 
$\yks= (\tau^{X,A}_n)^{-1}(\yk)$ and $\yks= \tau^{X,A}_n(\yk)$ respectively. 
Because of \eqref{tauAsup} we have  $\tau^{-1}B = (\tau^A_n)^{-1}B$
and $\tau B = \tau^A_n B$, so that after the substitution the 
indicator functions simplify to $1_B(\yk)$.  
Because of \eqref{gleichfallskongruentsup},
$\tau^{X,A}_n$ is constant on particles joined by bonds in $A$, 
so in $I_-$ we have
\begin{displaymath}
\begin{split}
\sigma_{x_1}&(\yk)  - \sigma_{x_2}(\yk) \; = 
\sigma_{x_1}(\tau_n^{X,A}(\yks)) - \sigma_{x_2}(\tau_n^{X,A}(\yks)) \\ 
&= \; \sigma_{x_1}(\yks) + \tau^{X,A}_n(x_1) 
  - \sigma_{x_2}(\yks) - \tau^{X,A}_n(x_2) \;  
= \; \sigma_{x_1}(\yks)  - \sigma_{x_2}(\yks),
\end{split}
\end{displaymath}
for every $\sigma \in (S^1)^{X}$ and for every bond $x_1x_2 \in
A$, and the same holds for $I_+$. 
Therefore the left hand side of \eqref{redsup} is equal to  
\begin{displaymath}
\begin{split}
&\int d\sigma_X \; \Big[ \; 1_{B}(\yk) \;  
 \prod_{x_1x_2 \in A} \; \left( \; 
 e^{J(x_1-x_2)v(\sigma_{x_1}(\yk) - \sigma_{x_2}(\yk))}-1 \; \right) 
 \; \\ 
&\hspace{1,5 cm} \cdot  \Big( \; \frac{e}{2}\;  
    e^{ -H^{\bar{U}}_{\Lambda_n}((\tau^{X,A}_n)^{-1}(\yk))} \; 
 + \; \frac{e}{2} e^{ -H^{\bar{U}}_{\Lambda_n}(\tau^{X,A}_n(\yk))} \;  
 - \; e^{ -H^{\bar{U}}_{\Lambda_n}(\yk)} \; \Big) \; \Big],
\end{split}
\end{displaymath}
which is nonnegative by \eqref{taylorsup} from  Lemma \ref{letaylorsup}. 
This proves \eqref{redsup} and completes the proof of the
 Lemma \ref{finishsup}.

%%%%%%%%%%%%%%%%%%%%%%%%%%%%%%%%%%%%%%%%%%%%%%%%%%%%%%%%%%%%%%%%
%%%%%%%%%%%%%%%%%%%%%%%%%%%%%%%%%%%%%%%%%%%%%%%%%%%%%%%%%%%%%%%%

\section{Proof of Theorem \ref{symallgsup}}

Let the assumptions of Theorem \ref{symallgsup} hold. 
First we observe that for every  $\tau \in G$ there is a torus $T$ 
such that $\tau \in T$ and $T$ is a subgroup of $G$. Every torus is a
finite product of compact 1-dimensional subgroups of $G$, so
w.l.o.g. we may assume that $\tau$ is contained in such a subgroup,
i. e. we may assume that $G$ is a compact 1-dimensional Lie-group, 
and hence that $G = S^1$ (for details see \cite{DS} for example). \\ 

For general $S$ we have to modify the decomposition of $V$. What we
need is a decomposition  $V = \bar{V} - v$ as guaranteed by Lemma
\ref{decompallgsup} presented below.\\

In order to deal with general $J$ we have to construct two different
decompositions of $V$: For $(x_1,\sigma_1),(x_2,\sigma_2) 
\in \R^2 \times S$ such that $J(x_1-x_2) \ge 0$ we decompose as
before: $V(\sigma_1,\sigma_2) = \bar{V}_{+}(\sigma_1,\sigma_2) 
- v_{+}(\sigma_1,\sigma_2)$, but if $J(x_1-x_2) < 0$ we decompose 
$V(\sigma_1,\sigma_2) = \bar{V}_{-}(\sigma_1,\sigma_2) 
+ v_{-}(\sigma_1,\sigma_2)$, where 
$v_{-}$ and $\bar{V}_{-}$ have the same properties as  $v_{+}$
and $\bar{V}_{+}$ respectively. This decomposition is also obtained
analogously to the following lemma.\\

\noindent The rest of the proof simply carries over. \qed\\
  
\noindent We still need

\begin{lemma} \label{decompallgsup}
Let $E$ be a compact topological space and let $S^1$ operate on $E$
continuously. Let $V: E^2 \to \R$ be a continuous mapping. Then we
have a decomposition $V = \bar{V} - v$ such that $ 0 < v <
\epsilon$, $\bar{V}$ is symmetric and $S^1$-invariant and such that
$\bar{V}(a,\tau b )$ is twice continuously differentiable with respect to
$\tau$ such that $\partial^2_{\tau} \bar{V}(a,\tau b)$ is bounded
uniformly in $a$ and $b$.
\end{lemma}

\Bew\\
Here we consider  $S^1 = \R / \Z$ and we identify functions on $S^1$ 
with periodic functions on $\R$. 
As a function of all three arguments $V(a,\tau b)$ is continuous on 
the compact space $E^2 \times S^1$, and therefore uniformly continuous.
Hence there exists a  $\delta > 0$ such that 
\begin{equation}\label{glmstetsup}
\alle a,b \in E \alle \tau',\tau \in \R: \quad 
|\tau' - \tau| < 2 \delta \; \Rightarrow \; 
|V(a,\tau' b) - V(a,\tau b)| < \frac{\epsilon}{2}.
\end{equation}
For this $\delta$  we choose a twice continuously
differentiable symmetric probability density $f_{\delta}: \R \to \R_+$ with
support in $[-\delta,\delta]$, for example  
\begin{equation*}
f_{\delta} (t) \; := \; c \cdot 1_{]-\delta,\delta[}(t) \cdot
 e^{-\frac{\delta^2}{\delta^2-t^2}} \quad \text{ with } \quad 
c \; := \; \int_{-\delta}^{\delta} e^{-\frac{\delta^2}{\delta^2-t^2}}\; dt .
\end{equation*}
Setting
\begin{equation*}
\bar{V}(a,b) \; := \; \int dt \;  f_{\delta}(t)  V(a,tb) + \frac{\epsilon}{2} 
\quad \text{ and } v \; := \; \bar{V} - V
\end{equation*}
gives us the desired decomposition.
$\bar{V}$ is measurable by Fubini's theorem, and symmetric, 
because $V$ is symmetric and $S^1$-invariant and
$f_{\delta}$ is symmetric. $\bar{V}$ is $S^1$-invariant because $V$ is. 
$0 < v < \epsilon$ is a straightforward consequence of \eqref{glmstetsup} 
and the small support of $f_{\delta}$. Finally, 
$\bar{V}(a,\tau b)  \; = \; \int dt \;  f_{\delta} (t-\tau)  V(a,tb) 
+ \frac{\epsilon}{2}$, which is twice continuously differentiable with respect
to $\tau$ such that $\partial^2_{\tau} \bar{V}(a,\tau b) \; 
= \; \int dt \;  f_{\delta}'' (t-\tau)  V(a,tb)$ 
is bounded by $2\delta \|f_{\delta}''\| \|V\|$. \qed 

\end{sloppypar}

%%%%%%%%%%%%%%%%%%%%%%%%%%%%%%%%%%%%%%%%%%%%%%%%%%%%%%
%%%%%%%%%%%%%%%%%%%%%%%%%%%%%%%%%%%%%%%%%%%%%%%%%%%%%%

\renewcommand{\thesection}{}

\setlength{\parindent}{0cm}

\end{document}